\documentclass[12pt,amstex]{article}
\usepackage[usenames]{color}
\usepackage{pgf,tikz}
\usetikzlibrary{arrows}

\usepackage{amscd}
\usepackage{amsmath}
\usepackage{amssymb}
\usepackage{graphicx}

\newtheorem {theo} {\bf Theorem} [section]
\newtheorem {prop} [theo] {\bf Proposition}
\newtheorem {cory} [theo] {\bf Corollary}
\newtheorem {lem} [theo] {\bf Lemma}
\newtheorem {defn} [theo] {\bf Definition}

\newtheorem {rem} [theo] {\bf Remark}
\newcommand{\QED}{\hfill \CaixaPreta \vspace{6mm}}
\def\CaixaPreta{\vrule Depth0pt height6pt width6pt}

\newcommand{\qed}{\nopagebreak\hfill{\vrule width6pt height6pt depth0pt}}

\newcommand{\dpy}{\displaystyle}
\newcommand{\be}{\begin{eqnarray}}
\newcommand{\ee}{\end{eqnarray}}
\newcommand{\benn}{\begin{eqnarray*}}
\newcommand{\eenn}{\end{eqnarray*}}
\newcommand{\bse}{\begin{equation}}
\newcommand{\ese}{\end{equation}}
\newcommand{\bsenn}{\begin{displaymath}}
\newcommand{\esenn}{\end{displaymath}}
\newcommand{\logand}{\;\;{\rm and }\;\;}

\newcommand{\logor}{\;\;{\rm or }\;\;}

\DeclareMathOperator{\Arg}{Arg}

\usepackage [height=9.5in,width=7in]{geometry}
\usepackage{graphicx}

\begin{document}

\title{Spectra of Tridiagonal Matrices}
\author{J. J. P. Veerman\thanks{Fariborz Maseeh Dept. of Math. and Stat., Portland State Univ.;
e-mail: veerman@pdx.edu} ,
D. K. Hammond\thanks{Applied Mathematics Department, Oregon Institute of Technology-Wilsonville; email: david.hammond@oit.edu} ,
 Pablo E. Baldivieso\thanks{Oregon State University - Cascades;
e-mail: pablo.baldivieso@osucascades.edu}}\maketitle

\vskip .0in

\noindent

\begin{abstract}
We characterize the eigenvalues and eigenvectors of a class of complex valued
tridiagonal $n$ by $n$ matrices subject to arbitrary boundary
conditions, i.e. with arbitrary elements on the first and last rows of
the matrix. 
For large $n$, we show there are up to $4$ eigenvalues, the so-called \emph{special eigenvalues},
whose behavior depends sensitively on the boundary conditions. The other eigenvalues,
the so-called \emph{regular eigenvalues} vary very little as function of the boundary
conditions. For large $n$, we determine the regular eigenvalues up to ${\cal O}(n^{-2})$, and the
special eigenvalues up to ${\cal O}(\kappa^n)$, for some $\kappa\in (0,1)$.
The components of the eigenvectors are determined up to ${\cal O}(n^{-1})$.

The matrices we study have important applications throughout the sciences.
Among the most common ones are arrays of linear dynamical systems with nearest neighbor
coupling, and discretizations of second order linear partial differential equations.
In both cases, we give examples where specific choices of boundary conditions substantially
influence leading eigenvalues, and therefore the global dynamics of the system.
\end{abstract}

\vskip 0.2in\noindent
\begin{centering}
\section{Introduction}
\label{chap:intro}
\end{centering}
\setcounter{figure}{0} \setcounter{equation}{0}

\noindent
We consider $n+1$ by $n+1$ complex valued tridiagonal
matrices given by the form
\begin{equation}
{A}_{n+1} =
\begin{pmatrix}
-{b}_0 & 1-{b}_1 & 0 & &\hdots &  0 \\
1 & 0 & 1 & &\hdots &   0 \\
0    & 1 & 0 & 1 & \hdots & 0 \\
\vdots &  & \ddots & \ddots& & \vdots\\
0 & \hdots & & 1 & 0 & 1 \\
0 & \hdots & & 0 & 1-{c}_{-1} & -{c}_0 \\
\end{pmatrix} \;.
\label{eq:Adef}
\end{equation}
In this paper, we characterize the spectrum and eigenvectors of $A$
when $n$ is large. Matrices of this type arise naturally when
describing the dynamics of systems of objects arranged in a line with
nearest-neighbor interactions, in this case the values of the
parameters $b_0,b_1,c_{-1}$ and $c_0$ are determined by how the
boundary conditions for the interactions are specified. A fundamental
question motivating this work is to understand how choices for the
boundary conditions can affect the global dynamics of such systems.

The results that we will obtain can be significantly extended by
allowing a few simple transforms. Suppose $D$ is a diagonal matrix
with non-zero diagonal entries in $\mathbb{C}$, $I$ the identity, and
$k\neq 0$ and $d$ are complex numbers.  Consider the following
transform of $A$:
\begin{equation}
B=k\left(DAD^{-1}+dI\right)
\label{eq:transformA}
\end{equation}
The spectrum of $B$ (see equation (\ref{eq:Bdef}) ) can be characterized in terms of the spectrum of $A$.
Details are given in Appendix 1.

Matrices of this form are commonly encountered in such a wide variety
of contexts that it is impossible to do the subject justice with a few
remarks.  They are found for example in one-dimensional arrays of coupled linear
ODE's whenever interactions are between nearest neighbors only. They
also occur in discretizations --- such as finite differences --- of
second order PDE's \cite{Haberman}. They are also important in solid
state physics where they play a crucial role in the study of crystal
vibrations (\cite{Ashcroft}, chapter 22). We will briefly discuss both
these examples in Section \ref{chap:applications}.  Many other uses
can be listed here. Our own interest derives from its uses in the
description of flocking and traffic (see for example \cite{Cantos2}).
We note that in many classical physics problems, the matrix $A$ must be symmetric.
In these cases, the boundary conditions are of course expressed in the first \emph{two} and
last \emph{two} rows of the matrix. This article does not cover those cases,
but see \cite{Fonsecaea} and references therein.

In \cite{Yueh} and \cite{Fonseca} special cases of the matrices
defined in equation (\ref{eq:transformA}) were studied. Eigenvalues of
tridiagonal matrices with the upper left block having constant values
were studied in \cite{Kulkarni1999}; this structure holds for our
matrix $A$ if $b_0=b_1=0$.  In \cite{Hammond}, the emphasis was on a
more general class of systems. There it was assumed that $b_1=0$, that
all parameters are real, and there were additional restrictions on the
diagonal matrix $D$ in equation (\ref{eq:transformA}). Here we
consider the general case where all parameters considered are
arbitrary complex numbers. There are also no restrictions on $D$
except that it must be invertible.

The structure of this paper is as follows. In section
\ref{chap:associated} we define a polynomial $H$ \emph{associated} to
the matrix $A$. This polynomial is not the characteristic
polynomial, but does have the property that the eigenvalues
$\lambda$ of $A$ are simple functions of the roots $r$ of $H$, namely
$\lambda=r+r^{-1}$.

In sections \ref{chap:regular} and \ref{chap:special} we give
approximate expressions for the roots of the associated
polynomial. The roots fall into two groups.  The ones we call
\emph{regular} (section \ref{chap:regular}) tend to fall close to the
unit circle (within $O(n^{-1})$). We determine them using topological
argument (Brouwer's fixed theorem) and we give expressions for them
that are accurate within $O(n^{-2})$.  In the other section we look at
the \emph{special} ones that fall ``far" from the unit circle and we
give expressions that are exponentially (in $n$) accurate. In section
\ref{chap:spectrum} we formulate our main theorem that gives the
eigenvalues of $A$.  In section \ref{chap:numerics} we describe
accurate numerical computation of eigenvalues based on these results,
and analyze the computational complexity.  Finally, section
\ref{chap:applications} discusses applications of these ideas to the
common physical assumption of periodic boundary conditions, and the
study of the eigenvalues of the discretized advection-diffusion
equation. In both of these applications, we show that for certain
parameter regimes the eigenvalue with largest real part (which is
necessarily significant for the global dynamics of the system) can be one of
the \emph{special} eigenvalues that strongly depends upon the boundary
conditions.

\vskip 0.2in\noindent
{\bf Acknowledgement:} We are grateful to Jeff Ovall for pointing out the usefulness
of the conjugation by a diagonal matrix (see Appendix 1).

\vskip 0.2in\noindent
\begin{centering}
\section{The Associated Polynomial}
\label{chap:associated}
\end{centering}
\setcounter{figure}{0} \setcounter{equation}{0}

\noindent
\begin{defn}
We define the $2n+4$ degree polynomial $H$ \underline{associated to $A$} as
\begin{equation*}
    H(z)=    z^{2n}(b_1 +{b}_0z + z^2)({c}_{-1} + {c}_0 z + z^2)
    - ({b}_1 z^2 + {b}_0 z + 1)({c}_{-1} z^2 + {c}_0 z + 1) ,
    \end{equation*}
and the \underline{auxiliary} functions $f$ and $g$ as
\begin{equation*}
f(z)= z^{2n} \quad \logand \quad
g(z)= \frac{({b}_1 z^2 + {b}_0 z +1)({c}_{-1}z^2 + {c}_0 z + 1)}
{({b}_1+{b}_0 z + z^2)({c}_{-1}+{c}_0 z + z^2)} .
\end{equation*}
Finally we define the \underline{auxiliary} polynomial
$p(z)=({b}_1 z^2 + {b}_0 z +1)({c}_{-1}z^2 + {c}_0 z + 1)$ and note that
\begin{equation*}
g(z)=\frac{p(z)}{z^4p(z^{-1})} \; .
\end{equation*}
\label{def:associated}
\end{defn}
We now describe how the eigenvalues of ${A}$ can be calculated by
analyzing the roots of $H$. In the following we denote the spectrum of
$A$ by $\sigma(A)$.

\vskip 0.2in\noindent
{\bf Remark:} If $b_1=1$ we see by inspection of $A_{n+1}$ that $-b_0$ is an eigenvalue
and that the remaining eigenvalues are equal to those of $A_{n}$ but now with $b_0$ set to 0
and $b_1$ to 0. We are thus allowed to assume without loss of generality that $b_1\neq 1$.
A similar remark holds for $c_{-1}$.

\vskip 0.2in\noindent
{\bf Remark:} The set of roots of $H$ is invariant under $z\rightarrow z^{-1}$.

\begin{prop}
Let $T$ be the set of roots of $\frac{H(y)}{(y-1)(y+1)}$. Then $\sigma(A)=\{ (y+y^{-1}) : y \in T\}$.
\label{prop:y2neq_deriv}
\end{prop}

\noindent
{\bf Proof:}
By the previous remarks we assume without loss of generality that $b_1\neq 1$ and $c_{-1}\neq 1$.

Note first that as $H(1)=H(-1)=0$, it follows that
$\frac{H(y)}{(y-1)(y+1)}$ is a polynomial. Letting $v=(v_0, v_1, ... , v_n)^T$, the eigenvalue equation
${A}v=rv$ is equivalent to the $n+1$ equations
\begin{align}
 &(1-{b}_1)v_1  = (r+{b}_0) v_0 \label{eq:bc1}\\
 & v_{k-1} + v_{k+1}  = rv_k
\mbox{ \quad for $1\leq k \leq n-1$} \label{eq:recurrence}\\
 &(1 - {c}_{-1}) v_{n-1}  = (r+{c}_0) v_n \label{eq:bc2}
 \end{align}

 We will proceed by writing the general solution to the linear
 recurrence relation implied by (\ref{eq:recurrence}), with the other
 two equations above providing boundary conditions.
 Equation
(\ref{eq:recurrence}) implies $v_{k+1} = r v_k- v_{k-1}$.
 Introducing $C=\begin{pmatrix} 0 & 1 \\ -1 & r \end{pmatrix}$, we may
 rewrite this as
 $\begin{pmatrix}v_k \\ v_{k+1}\end{pmatrix}=C\begin{pmatrix} v_{k-1}\\
   v_k\end{pmatrix}$, which implies $\begin{pmatrix}v_k \\
   v_{k+1}\end{pmatrix} = C^k \begin{pmatrix} v_0 \\
   v_1 \end{pmatrix}$.

 The characteristic polynomial of $C$ is $\lambda^2 - r \lambda +1$,
 which has repeated roots precisely when $r^2-4=0$. Thus if
 $r\neq \pm2$, the eigenvalues of $C$ will be distinct. Assume for now
 this is the case, and denote the eigenvalues of $C$ by $x_+$ and
 $x_-$. It then follows that we must have
 $v_k = c_+ x_+^k + c_- x_-^k$ for $0\leq k \leq n$, for some
 constants $c_+$ and $c_-$. Valid eigenvalues $r$ will be those such
 that these expressions for $v_k$ are also consistent with the
 boundary conditions (\ref{eq:bc1}) and (\ref{eq:bc2}). 
 These imply
\begin{align*}
(1-{b}_1)(c_+x_+ + c_-x_-) &= (r+{b}_0) (c_+ + c_-) \\
(1-{c}_{-1})(c_+x_+^{n-1} + c_-x_-^{n-1}) &= (r+{c}_0) (c_+x_+^n +
                                            c_-x_-^n) \;.
\end{align*}

We now note that $x_+ + x_- = \textrm{trace}(C) = r$ and
$x_+ x_- = \det(C) = 1$. We use the latter to introduce the
substitution $x_+= y$ and $x_-=y^{-1}$. These then imply that
\begin{equation}
r=(y+y^{-1}) \quad \logand \quad v_k = c_+ x_+^k + c_- x_-^k \; .
\label{eq:eigenpair}
\end{equation}
Substituting these into the above and simplifying gives the system of equations
(noting that $y\neq 0$)
\begin{equation}
\begin{pmatrix}
{b}_1y^2 + {b}_0 y + 1 & {b}_1 + {b}_0 y + y^2\\
y^{2n}({c}_{-1} +{c}_0y + y^2) & c_{-1}y^2 + {c}_0y  + 1\\
\end{pmatrix}
\begin{pmatrix} c_+ \\ c_-\\ \end{pmatrix} =
\begin{pmatrix} 0 \\ 0\\ \end{pmatrix} \;.
\label{eq:cplus-cminus}
\end{equation}
There will be nontrivial solutions for $c_\pm$ if and only if the determinant of the
corresponding matrix is zero.
This corresponds exactly to $H(y)=0$. If $y$ is a root of the above equation
not equal to $\pm 1$, then the corresponding $r=y+y^{-1} \neq \pm 2$,
and the previous steps imply that $r\in \sigma(A)$.

We now consider the case when $r=\pm 2$. We show that this occurs
exactly when $H(y)$ has a repeated root at $\pm 1$, so that the polynomial
$h(y)\equiv \frac{H(y)}{(1-y)(1+y)}$ will have a root at $y=\pm 1$.
Denote $\xi=(b_0,b_1,c_0,c_{-1}) \in \mathbb{C}^4$. Clearly $h_\xi(z)$ is a polynomial
evaluated in $z$ and so is a continuous function of $\xi$.
It is also well-known that the eigenvalues of a matrix are continuous functions of its entries
(in this case $\xi$). Let $y_+(\xi)$ be $\frac r2 +\sqrt{\frac{r^2}{4}-1}$. Choose a path
$\xi(t)$ so that $r(\xi(t))=\pm2$ iff $t=0$. By continuity we have
\begin{equation*}
\lim_{t\rightarrow
  0}\,h_{\xi(t)}(y_+(\xi(t)))=h_{\xi(0)}(y_+(\xi(0)))=0 \; ,
\end{equation*}
and thus the polynomial $h_{\xi(0)}$ has a root at $\pm1$.
\QED

\vskip .0in\noindent
\begin{centering}
\section{Regular Roots of the Associated Polynomial}
\label{chap:regular}
\end{centering}
\setcounter{figure}{0} \setcounter{equation}{0}

\noindent
We begin our study of the roots of $H(y)$ of Definition \ref{def:associated}.
First we introduce the following notation.

\vskip -0.0in\noindent
\begin{defn} Let $\gamma(t)= e^{it}$ for $t\in[0,2\pi)$. Using the auxiliary
function $g$ from Definition \ref{def:associated}, we define differentiable functions $R$
and $\Psi$ from $\mathbb{R}$ to $\mathbb{R}$ by requiring:
\begin{equation}
g( e^{it})=R(t)e^{i\Psi(t)}
\label{eq:RandPsi}
\end{equation}
Assume that $g$ has no zeros or poles on the unit circle.
Then $e^{i\Psi(t)}$ and $g\big|_\gamma$ have the same well-defined winding number $w\in \mathbb{Z}$.
\label{def:RandPsi}
\end{defn}

\noindent
{\bf Remark:} Note that the continuous map $\Psi$ is the \emph{lift} of
$e^{i\Psi}:\mathbb{R}\rightarrow S^1$ (the circle)
to the real line. Thus $\Psi(2\pi)-\Psi(0)=2 \pi w$.


\begin{defn} Let $Q$ be the number of zeros (with multiplicity) of the
  auxiliary polynomial $p$ inside the unit circle.
\label{def:Q}
\end{defn}

\begin{lem}
The winding number $w$ of $g\big|_\gamma$ equals 2Q-4.
\label{lem:windingnr}
\end{lem}

\noindent
{\bf Proof:} The winding number satisfies (see \cite{apostol}):
\begin{equation*}
w=\frac{1}{2\pi i}\,\int_\gamma\, \frac{g'(z)}{g(z)}\, dz=N-P
\end{equation*}
where $N$ is the number of zeroes of $g$ inside $\gamma$ and $P$ the number of poles inside $\gamma$.
Clearly $N=Q$. Furthermore $z^4p(z^{-1})$ is a quartic polynomial whose roots are the inverses
of the roots of $p(z)$. Hence $P=4-Q$. \QED

\begin{defn}
Choose $\Delta$ such that $\Delta^{-1}<1<\Delta$ and let
$A_\Delta=\{z\in \mathbb{C}\,|\,\Delta <|z|<\Delta\}$.
Choose $C>1$ and $D>0$ be constants so that on $A_\Delta$:
\begin{equation*}
C^{-1} < | g(z)| < C \quad \logand \quad   |g'(z)| < D \; .
\end{equation*}
\label{def:annulus_bounds}
\end{defn}

\vskip -0.3in
We now collect a few lemmas. The first two proofs are elementary and are left to the reader.

\begin{lem} For all $x>0$: $\quad \dpy 1-x^{-\frac{1}{2n}}\leq x^{\frac{1}{2n}}-1
\leq \dfrac{x-1}{2n}\quad $ (with equality iff $x=1$).
\label{lem:calculus-inequality}
\end{lem}

\begin{lem} Let $f$ the auxiliary function of Definition \ref{def:associated}, then on $A_\Delta$ we have:
\begin{equation*}
\frac{1}{2n(1+\Delta)} \leq
\left|\frac{d}{dz}f^{-1}(z)\right| \leq \frac{1}{2n(1-\Delta)}.
\end{equation*}
\label{lem:f-inverse-prime}
\end{lem}

\begin{lem} The function $\Psi$ satisfies $\quad \dpy \left|\Psi'(t)\right|^2 =
\dfrac{\left|g'(e^{it})\right|^2-\left|R'(t)\right|^2}{\left|g(e^{it})\right|^2} < C^2D^2 \quad $.
\label{lem:Psi-prime}
\end{lem}

\noindent
{\bf Proof:} Differentiating equation (\ref{eq:RandPsi}) with respect to $t$ gives
\begin{equation*}
g'\left( e^{it}\right) e^{it}\, i=\left[R'(t)+i\Psi'(t)R(t)\right]\,e^{i\Psi(t)}.
\end{equation*}
Taking the absolute value and squaring yields
\begin{equation*}
\left| g'\left( e^{it}\right)\right|^2=R'(t)^2+\Psi'(t)^2R(t)^2,
\end{equation*}
and with $R(t)=|g(e^{it})|$ and Definition \ref{def:annulus_bounds} this implies the result.
\QED

\begin{defn} We call $t\in[0,2\pi)$ a \underline{phase root} if it is a
  solution to $e^{2int}=e^{i\Psi(t)}$.
\label{def:phase_roots}
\end{defn}

\begin{prop} \label{prop:the_phase_roots}
  For any value of $n$, there are at least $2n+4-2Q$ phase roots in
  $[0,2\pi)$. Furthermore, for $n>\frac{CD}{2}$, there are exactly
  $2n+4-2Q$ phase roots $\{t_i\}_{i=1}^{2n+4-2Q}$ in $[0,2\pi)$, and
  these roots ordered in ascending magnitude satisfy
\begin{equation*}
\dfrac{2\pi}{n+CD}\leq |t_k-t_{k-1}|\leq \dfrac{2\pi}{n-CD} .
\end{equation*}
\label{prop:phase-roots}
\end{prop}

\noindent
{\bf Proof:}
Let $h$ be the continuous function from $\mathbb{R}$ to itself given
by $h:t\rightarrow 2nt-\Psi(t)$. The phase roots are the solutions of
$h(t)\in 2\pi \mathbb{Z}$.  Since $h$ satisfies
$h(2\pi)-h(0)=2\pi(2n-w)$ (see the remark after Definition
\ref{def:RandPsi}), we have at least $2n-w$ solutions. This is equal
to $2n+4-2Q$ by Lemma \ref{lem:windingnr}, which proves the first
assertion of the proposition.

Now assuming the condition on $n$, Lemma \ref{lem:Psi-prime} gives
\begin{equation*}
0<2n-CD<h'(t)<2n+CD,
\end{equation*}
and thus there are exactly $2n-w$ solutions. Also for two successive phase roots $t_{k+1}$ and
$t_k$ the mean value theorem gives
\begin{equation*}
2n-CD<\dfrac{2\pi}{t_{k+1}-t_k}<2n+CD,
\end{equation*}
which implies the result.
\QED

Before we present our main result, which is valid for complex
matrices, we address an important remark for the real matrix case.

\begin{prop}\label{prop:phase_roots_real}
If $A$ is real, then each $e^{it_k}$ is an exact
root of $H(z)$, where $t_k$ is a phase root.
\end{prop}

\noindent
{\bf Proof:}
If $A$ is real-valued then the coefficients of $p(z)$ are real, and so for any
$t\in \mathbb{R}$ one has $\overline{p(e^{it})}=p(e^{-it})$. It
follows that $|p(e^{it})|=|p(e^{-it})|$, which implies $|g(e^{it})|=1$.
Now let $t_k$ be a phase root. Then $e^{i 2 n t_k} = e^{i\Psi(t_k)}$
where $\Psi$ is the lift map defined in Equation (\ref{eq:RandPsi}). But as
$|g(e^{it})|=1$ for all $t$, it must be that
$e^{i\Psi(t_k)} =g(e^{it_k})$ and thus that $e^{it_k}$ is an exact root
of $H(z)$. \QED

\begin{defn} For $g$ as in Definition \ref{def:associated} and each
  $1\leq k \leq 2n+4-2Q$, define the \underline{approximate roots}
\begin{equation*}
z^*_k =|g( e^{i t_k})|^{\tfrac{1}{2n}} e^{i t_k}
\end{equation*}
and the \underline{approximation discs}
\begin{equation*}
  B_k = \left\{ z\, : \,|z-z_k^*| \leq  \frac{CD}{2(1-\Delta)n^2}\right\}.
\end{equation*}
\label{def:approx-discs}
\end{defn}

\begin{theo} Suppose $n>\max\left\{\frac{CD}{2(1-\Delta)}, \frac{C}{\Delta}\right\}$.
Then each of the $2n+4-2Q$ discs $B_k$
contains a unique root of $H(z)$ in Definition \ref{def:associated}.
\label{thm:regular-roots}
\end{theo}

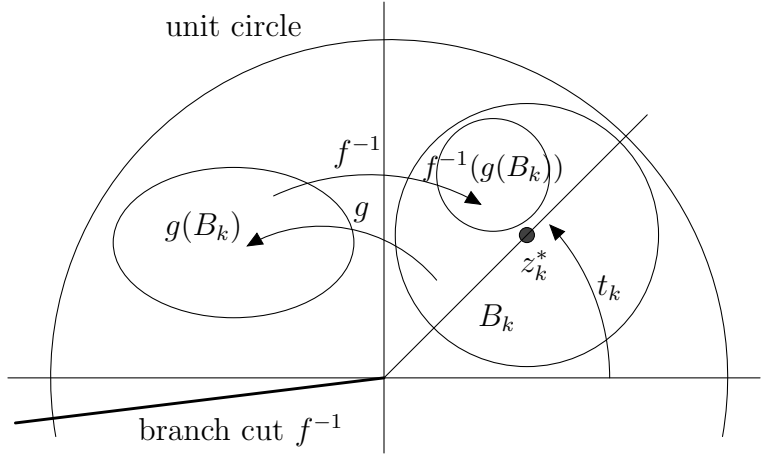
\begin{figure}[pbth]
\center
\begin{tikzpicture}[line cap=round,line join=round,>=triangle 45,x=1cm,y=1cm]
\draw[color=black] (-5,0) -- (5,0);
\draw[color = black] (0,5) -- (0,-1);
\draw[line width=1.1pt,fill=black,fill opacity=0.75] (0,0) -- (-4.9,-.6);

\draw (4.5,-.781) arc (-10:190:4.5cm);
\draw(1.45,2.7) circle (.75cm);

\draw[color=black] (0,0) -- (3.5,3.5);

\draw(1.9,1.9) circle (1.75cm);
\draw [fill=black,fill opacity=0.75] (1.9,1.9) circle (0.1cm);

\draw (2,1.5) node {$z^*_k$};
\draw (3.,1.2) node {$t_k$};
\draw (1.5,.8) node {$B_k$};
\draw (-2.4,2) node {$g(B_k)$};
\draw (-.3,2.2) node {$g$};
\draw (-.35,3) node {$f^{-1}$};
\draw (1.45,2.8) node {\small{$f^{-1}(g (B_k))$}};
\draw (-1.9,-.7) node {branch cut $f^{-1}$};
\draw (-2,4.7) node {unit circle};

\draw (-2,1.8) ellipse (1.6cm and 1.cm);

\draw[->] (0:3.cm) arc (0: 43:3.cm);

\draw[->] (.7,1.3) arc (40: 120:2.cm);

\draw[<-] (1.3,2.3) arc (60: 115:3.cm);

\end{tikzpicture}

\caption{ \emph{ Illustration of the phase roots $t_k$, approximate roots
  $z_k^*$ and approximation discs $B_k$ from the proof of Theorem \ref{thm:regular-roots}. }}
\label{fig:brouwer}
\end{figure}

\noindent {\bf Proof:} Let $f$ and $g$ be as in Definition
\ref{def:associated}. We show that for each $k\in \{1,\cdots,2n-w\}$
we can choose an inverse branch $f^{-1}$ of $f$ such that
$f^{-1}\circ g$ is continuous and maps B$_k$ to B$_k$ (see Figure
\ref{fig:brouwer}). By Brouwer's theorem this gives a fixed
point. Uniqueness is then implied by the observation that on $B_k$,
$f^{-1}\circ g$ is a contraction.
Let $z\in B_k$, then
\begin{equation*}
|z-e^{it_k}|\leq |z-z_k^*| + |z_k^*-e^{it_k}| \leq \frac{CD}{2(1-\Delta)n^2} +
\left(|g( e^{it_k})|^\frac{1}{2n}-1\right) < \frac{CD}{2(1-\Delta)n^2}+\frac{C-1}{2n}.
\end{equation*}
For the last inequality we have used Definition \ref{def:approx-discs}
and Lemma \ref{lem:calculus-inequality}. By the hypothesis on $n$,
this last quantity is less than $\frac Cn$ which in turn is less than
$\Delta$ and thus $B_k \subseteq A_\Delta$.

We can thus use Definition \ref{def:annulus_bounds} to ensure that for $z\in B_k$
\begin{equation*}
  |g(z)-g(e^{it_k})| < \frac{CD}{n} \;.
\end{equation*}
By Definition \ref{def:RandPsi} and Proposition \ref{prop:phase-roots} we have that
$z_k^*=f^{-1}\circ g( e^{it_k})$. With the above equation and using Lemma
\ref{lem:f-inverse-prime} this gives
\begin{equation}
\big| f^{-1}\circ g(z)-z_k^*\big|=\big| f^{-1}\circ g(z)-f^{-1}\circ g( e^{it_k})\big|
<\frac{DC}{2(1-\Delta)n^2},
\label{eq:contraction}
\end{equation}
which proves that $f^{-1}\circ g(B_k)\subseteq B_k$.

Since $B_k\subseteq A_\Delta$ we have that on $B_k$, $|g(z)|>C^{-1}$.
Thus if $g(B_k)$ encircles the origin, there must be points $z_1$ and $z_2$ in $B_k$ such that
by using Definition \ref{def:annulus_bounds}
\begin{equation*}
2C^{-1}< |g(z_2)-g(z_2)|\leq D \operatorname{diam}(B_k)=\frac{CD^2}{2(1-\Delta)n^2}.
\end{equation*}
But this is impossible by hypothesis. Thus we choose a branch cut for
$f^{-1}$ so that the local inverse on $g(B_k)$ is continuous. This
establishes the existence of the fixed point.

The fact that $f^{-1}\circ g$ is a contraction on $B_k$ follows from this simple calculation:
\begin{equation} \label{eq:contraction_map_derivative}
\Big|\frac{d}{dz}\,f^{-1}\circ g(z)\Big|=\Big|\frac{d}{dz}\,f^{-1}\Big|_{g(z)}\Big|\cdot
|g'(z)|< \frac{D}{2n(1-\Delta)},
\end{equation}
which is smaller than 1 by the hypothesis on $n$.  \QED

\noindent {\bf Remark:} The contraction mapping $f^{-1}\circ g$ can be
iterated to give more accurate estimates of the roots, this is
developed further in section \ref{chap:numerics}.

\vskip .2in\noindent
\begin{centering}
\section{Special Roots of the Associated Polynomial}
\label{chap:special}
\end{centering}
\setcounter{figure}{0} \setcounter{equation}{0}

\noindent
If $n$ is sufficiently large so that Theorem \ref{thm:regular-roots} holds,
then we refer to the $2n+4-2Q$ roots of $H(z)$ that are contained in the
approximation discs $B_k$ as ``regular roots'', the remaining roots of
$H(z)$ will be denoted as ``special roots''.

\vskip .2in\noindent
\begin{prop}
Let $z_0$ be a root of $g(z)$ that is inside the unit circle, with
multiplicity $m$, and fix $\rho$ satisfying $|z_0| < \rho < 1$.  Then
there is a constant $K$ such that for sufficiently large $n$ the
circle of radius $\epsilon =K (\rho^{1/m})^{2n}$ centered at $z_0$ contains
$m$ roots of roots of $z^{2n}-g(z)$.
\label{prop:rouche}
\end{prop}

\noindent
{\bf Proof:}
We apply Rouch\'{e}'s theorem (see \cite{ahlfors}) to $f_1(z) = g(z)$ and $f_2(z) = g(z)-z^{2n}$.

Pick an $\epsilon$ so that $0<\epsilon<\rho-|{z_0}|$ and denote $D_{z_0}(\epsilon)$ the sphere of
radius $\epsilon$ centered at ${z_0}$. On $D_{z_0}(\epsilon)$ we have
\begin{eqnarray*}
&&|f_1(z)- f_2(z)| = |z|^{2n} < \rho^{2n}\\
&&|f_1(z)| = \Big| \frac{g^{(m)}({z_0})}{m!} (z-{z_0})^m + O((z-{z_0})^{m+1})|>M \epsilon^m
\end{eqnarray*}
for $M=\big| \frac 12 \frac{g^{(m)}({z_0})}{ m!}\big|$. Thus if we set $\epsilon = M^{-1/m} \rho^{2n/m}$,
we have that $|f_1(z)-f_2(z)|<|f_1(z)|$. Hence by Rouch\'{e}'s theorem $f_1(z)=g(z)$ and
$f_2(z)=z^{2n}-g(z)$ must have the same number of zeros in $D_{z_0}(\epsilon)$. \QED

\vskip .0in\noindent
\begin{theo} If p(z) has $Q$ roots inside the unit circle, then for $n$ large enough
$H(z)$ in Definition \ref{def:associated} has $2Q$ special roots (counting algebraic multiplicity).
\label{thm:special-roots}
\end{theo}

\noindent {\bf Proof:} Proposition \ref{prop:rouche} shows there are
$Q$ roots of $H(z)$ associated with the $Q$ roots of $g(z)$ inside the
unit circle. Since the set of roots is invariant under
$z\rightarrow z^{-1}$, there must also be $Q$ roots outside the unit
circle. For large enough $n$, none of these roots are in the
approximation discs $B_k$ of Definition \ref{def:approx-discs} since
these discs can be made to lie arbitrarily close to the unit circle as
$n\rightarrow \infty$. Finally we note that all roots of $g(z)$ are
roots of $p(z)$ (see Definition \ref{def:associated}).  \QED

\vskip .2in\noindent
\begin{centering}
\section{Eigenvalues and Eigenvectors of $A$}
\label{chap:spectrum}
\end{centering}
\setcounter{figure}{0} \setcounter{equation}{0}

\noindent
The main result below is an immediate corollary of Theorems \ref{thm:regular-roots} and
\ref{thm:special-roots}. By Proposition \ref{prop:y2neq_deriv} each of these two sets of roots
is invariant under $z\rightarrow z^{-1}$. In that same Proposition we see that 2 roots
$y$ and $y^{-1}$ combine to give an eigenvalue.

\vskip .2in\noindent
\begin{cory}
Let the parameters for the matrix $A$ be such that the auxiliary function $g(z)$ has no zeros or
poles on the unit circle. Then, for sufficiently large $n$, the spectrum of $A$ consists of
$n+1-Q$ \underline{regular} eigenvalues $\{r_k\}_{k=1}^{n+1-Q}$ and $Q$
\underline{special} eigenvalues $\{s_k\}_{k=1}^{Q}$, given by
\begin{align*}
r_k & = |g(e^{i t_k})|^{1/2n} e^{it_k}+ |g(e^{i t_k})|^{-1/2n}  e^{-it_k} + O(n^{-2})\\
s_k & = y_k + y_k^{-1} + O(\kappa^{-2n})
\end{align*}
where $t_k$ are the $n+1-Q$ phase roots satisfying $0<t_k < \pi$, and $y_k$ are the $Q$ roots of
the auxiliary polynomial $p(z)$ inside the unit circle, and $\kappa$ is a number greater than 1.
\label{cor:main}
\end{cory}

\vskip .1in\noindent
This result allows us to determine the eigenvectors of $A$. Let $z$ be one of the regular
roots of Theorem \ref{thm:regular-roots}, then equations (\ref{eq:eigenpair})
and (\ref{eq:cplus-cminus}) imply that the components $v_k$ of the eigenvector associated
to $y+y^{-1}$ are given by
\begin{equation}
v_k= (b_1+b_0z+z^2)\,z^k - (b_1z^2+b_0z+1)\,z^{-k} + O(n^{-1}) \;.
\label{eq:component}
\end{equation}
The error of $O(n^{-1})$ in $z^k$ for $k\in\{-n,\cdots n\}$ follows because $z$ itself
is determined up to $O(n^{-2})$. The modulus $|v_k|$ is bounded in some interval $[K^{-1},K]$
for some $K>1$ independent of $n$.

The eigenvectors associated with any special root $z$ exhibit a different behavior.
In this case the error in $z$ is exponentially small in $n$ (see Theorem \ref{thm:special-roots}).
Thus the error in  $z^k$ for $k\in\{-n,\cdots n\}$ is also exponential. Equation
(\ref{eq:component}) holds but with an error $O(\tau^{2n})$ for some $0<\tau<1$.
However in this case the values $|z^k|$ become exponentially large and those of
$|z^{-k}|$ exponentially small (or vice versa, depending on the value of $|z|$).

Finally in the case that we have an eigenvalue $\pm 2$, the eigenvectors of A are well-known:
we have
\begin{equation*}
 v_k=-(k-1)a+kb
\end{equation*}
for arbitrary $a$ and $b$. This of course is exact.

Our results simplify in the important case when $A$ is real valued. In
this case the ``approximate roots'' $z_k^*$ from Definition
\ref{def:approx-discs} are in fact exact.  Additionally, we may remove
the additional assumption that $g(z)$ has no roots or poles on the
unit circle.

\vskip .2in\noindent
\begin{cory}
Suppose the matrix $A$ is real. Then $g(z)$ has no 
zeros or poles on the unit circle, and
for sufficiently large $n$, the spectrum of $A$ consists of
$n+1-Q$ \underline{regular} eigenvalues $\{r_k\}_{k=1}^{n+1-Q}$ and $Q$
\underline{special} eigenvalues $\{s_k\}_{k=1}^{Q}$, given by
\begin{align*}
r_k & = 2\cos t_k  \\
s_k & = y_k + y_k^{-1} + O(\kappa^{-2n})
\end{align*}
where $t_k$ are the $n+1-Q$ phase roots satisfying $0<t_k < \pi$, and $y_k$ are the $Q$ roots of
$p(z)$ inside the unit circle, and $\kappa$ is a number greater than 1.
\label{cor:main-real}
\end{cory}

\vskip .0in\noindent
{\bf Proof:} If $A$ is real then the coefficients of $p(z)$
are real, so roots of $p(z)$ are either real or occur in conjugate
pairs. Thus if $p(z)$ has a root $e^{i\phi}$ on the unit circle, then
$e^{-i\phi}$ is also a root of $p$. However this implies that
$z^4 p(z^{-1})$ has a root at $e^{i\phi}$, and so these roots cancel
in the expression for $g(z)$. Thus $g(z)$ can have no zeros or poles
on the unit circle.
Next, Proposition \ref{prop:phase_roots_real} implies that the phase
roots $t_k$ yield exact roots $e^{it_k}$ of $H(z)$. The correponding
eigenvalues $r_k=e^{it_k} + 1/e^{it_k}$ (notably, without the
$O(n^{-2})$ term) imply the desired result.
\QED

\setcounter{figure}{0} \setcounter{equation}{0}

\vskip .0in\noindent
\begin{centering}
\section{Numerical Eigenvalue Computation }
\label{chap:numerics}
\end{centering}

We describe and analyze a numerical procedure for computing the
regular eigenvalues of $A$ to machine precision based on first
computing the phase roots then iterating the contraction mapping
described in section \ref{chap:regular}.  We compute the phase roots
by applying the bisection method to determine the roots of
$k(t)=\Arg(\frac{e^{i2nt}}{g(e^{it})})$, where
$\Arg : \mathbb{C}\setminus \{0\} \to (-\pi,\pi]$ gives the angle of a
complex number.
The bisection method is guaranteed to converge to a root, if
initialized with the endpoints of an interval (called a bracket) that
contains a root and over which $k(t)$ is continous. We determine a set
of brackets for the phase roots by setting $N=6n$ and defining the $N$
intervals $u_\ell=\frac{2\pi}{N} \ell$ for $1\leq \ell \leq N+1$.
Note that the function $k(t)$ is pointwise discontinuous at any value
$t$ where $k(t)=\pi$.  We retain as brackets the intervals
$I_\ell = [u_\ell, u_{\ell+1}]$ for which
$k(u_\ell) k(u_{\ell+1}) < 0$, and for which
$|k(u_\ell) - k(u_{\ell+1})|<\pi/2$. This latter condition is needed
to avoid retaining brackets which contain a point where $k(t)$ is
discontinuous.


\vskip .2in\noindent
\begin{prop} \label{prop:brackets}
  Let $n>CD/2$, where $C$ and $D$ are from definition
  \ref{def:annulus_bounds}. Then, for $N\geq 6n$, each interval $I_\ell$
  can contain at most one phase root. Additionally, each interval
  $I_\ell$ which contains a point of discontinuity of $k(t)$ will
  satisfy $|k(u_\ell) - k(u_{\ell+1})|>\pi/2$.
\end{prop}

\noindent
{\bf Proof:}
Proposition \ref{prop:the_phase_roots} implies that the distance
between any successive two phase roots is at least
$\frac{2\pi}{n+CD}$. As the length of $I_\ell$ is $\frac{2\pi}{N}$,
$I_\ell$ cannot contain two phase roots if
$\frac{2\pi}{N} <\frac{2\pi}{n+CD}$. This is ensured provided
$N>n+CD$, which is ensured by the assumption on $n$ for any $N>3n$.

Second, observe that $k(t) = h(t) (\mod 2\pi)$, where
$h(t)=2nt-\Psi(t)$ is as defined in the proof of
\ref{prop:the_phase_roots} and where representative angles are chosen
on $(-\pi,\pi]$. As $h'(t)<2n+CD$, over a single interval, $h(t)$ can
change by no more than
$(2n+CD)\frac{2\pi}{N} <  \frac{8\pi n}{N} < \frac{3\pi}{2}$, where the
last inequality follows from $N>6n$. This implies that if an interval
does contain a jump discontinuity (at which $k(t)$ changes by $2\pi$),
the values of $k$ at the endpoints will differ by more than $\pi/2$.
\QED

Given the  $2n+4-2Q$ phase roots $t_k$, we define $f_k^{-1}$ to be a
 branch\footnote{ Explicitly, we define
  $f_k^{-1}(re^{i\theta}) = r^{1/2n} e^{i\theta'/2n}$, where
  $\theta'=\theta+2\pi m$ and $m$ is such that
  $2nt_k - \pi \leq \theta' < 2nt_k + \pi$. This holds for
  $m=\lceil \frac{2nt_k-\theta-\pi}{2\pi}\rceil$.}
of the inverse of $f(z)=z^{2n}$ satisfying $f_k^{-1} ((e^{i t_k})^{2n}) =e^{it_k}$.
Define the iterates $z_k^{(i)}$ by setting
\begin{equation} \label{eq:contraction_iter}
z_k^{(i)} = f_k^{-1} \circ g (z_k^{(i-1)})
\end{equation}
with $z_k^{(0)} = e^{it_k}$. From the analysis in section
\ref{chap:regular}, it follows that
$\lim_{i\to\infty} z_k^{(i)}=\tilde{z}_k$ is a root of
$z^{2n}=g(z)$.

Given a fixed desired precision $\epsilon$, our
numerical procedure for computing the regular eigenvalues consists of
the steps: (1) Compute the phase roots $t_k$ to within $\epsilon$ by
bisection (2) For each phase root, iterate equation
(\ref{eq:contraction_iter}) until convergence within $\epsilon$ (3)
Calculate the eigenvalues via $r=z + 1/z$ where $z$ is the converged
result from step 2.

\begin{figure}
\centerline{
\includegraphics[width=7cm]{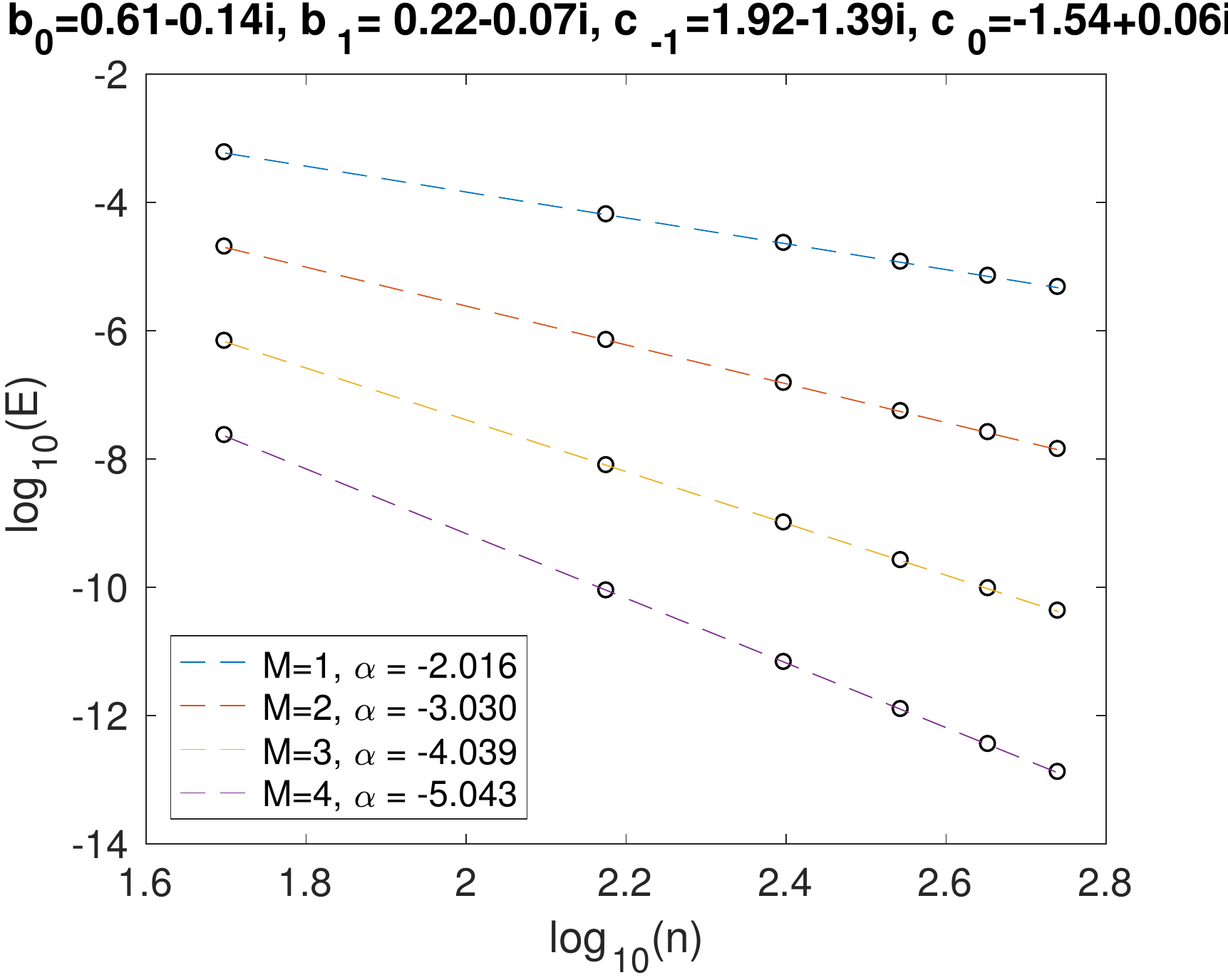} }
\caption{\emph{Maximum error in estimated eigenvalues. Errors are computed
  by comparing eigenvalues computed by the proposed method, using a fixed
  number $M$ of iterations of (\ref{eq:contraction_iter}), to those
  computed using the numerical eigenvalue routine eig in
  MATLAB. Given $\alpha$ values are the slopes of the least-squares
  linear fits (dotted lines).    }}
\label{fig:error_vs_n}
\end{figure}

Proposition \ref{prop:brackets} implies that, for $N=6n$ , this
procedure is guaranteed to find all of the regular eigenvalues of $A$
(provided $n>CD/2$). Before discussing the computational complexity of
this procedure, we analyze the iterates of equation
(\ref{eq:contraction_iter}).  Theorem \ref{thm:regular-roots} implies
that $|z^{(1)}_k-\tilde{z}_k|<\frac{CD}{2(1-\Delta)n^2}$.  Using the
bound on $(f^{-1}\circ g)'$ from equation
(\ref{eq:contraction_map_derivative}), we see that the later iterates
satisfy
\begin{equation}\label{eq:contraction_iter_bound}
|z_k^{(i)} - \tilde{z}_k| < 2C \left(\frac{D}{2(1-\Delta)}\right)^i
\frac{1}{n^{i+1}} .
\end{equation}
This implies that the residual error in the eigenvalues computed from
applying $M$ steps of (\ref{eq:contraction_iter}) is proportional to
$\frac{1}{n^{M+1}}$.
This behavior is illustrated in Figure \ref{fig:error_vs_n}, where we show the
maximum error of the estimated eigenvalues as computed by our method
vs $n$, for $n\in \{50,150,250,350,450,550\}$, and $M\in\{1,2,3,4\}$.
On a log-log plot, the observed slopes are close to $-(M+1)$,
consistent with error proportional to $\frac{1}{n^{M+1}}$.

We now examine the computational complexity of our overall numerical
procedure as a function of the matrix size $n$, by counting the number
evaluations of either $k(t)$ or $f^{-1} \circ g$ as a proxy for
computational cost.
Computation of the initial brackets takes $N=6n$ function
evaluations. The error from the bisection method after $q$ steps is
bounded by $2^{-q}$ times the length of the original bracketing
interval, in our case $\frac{2\pi}{6n}$. This implies the need for
$q=\log_2(\epsilon^{-1})-\log_2(3n)$ bisection steps for each phase
root, implying a total cost of $6n+(2n+4-2Q)
(\log_2(\epsilon^{-1})-\log_2(3n))$ to compute all of the phase roots.
Iterating equation (\ref{eq:contraction_iter}) for all of the roots
requires at most $(2n+4-2Q)M$ function evaluations, where $M$ is the
maximum number of iterations performed.  The bound
(\ref{eq:contraction_iter_bound}) implies that convergence within
$\epsilon$ is assured if
$\frac{2C}{n} (\frac{D}{2(1-\Delta)n})^M < \epsilon$.  Under the
conditions $n>2C$ and $n>\frac{D}{(4(1-\Delta)}$, which hold for
sufficiently large $n$, convergence within $\epsilon$ is assured for
$(\frac{1}{2})^M < \epsilon$, which holds for
$M=\log_2(\epsilon^{-1})$.  Thus for sufficiently large $n$, iterating
equation (\ref{eq:contraction_iter}) for all of the roots will require
no more than $\log_2(\epsilon^{-1}) (2n+4-2Q)$ function evaluations.
Together, these imply that the total computational cost of computing
all of the regular eigenvalues is bounded by
$6n+(2n+4-2Q)(2\log_2(\epsilon^{-1})$, which is
$O(n\log_2(\epsilon))$. For $\epsilon$ fixed independent of $n$,
the overall computational complexity of our approach is $O(n)$. This
should be contrasted with the standard QZ algorithm for computing all of the
eigenvalues of a matrix, which has complexity $O(n^3)$ \cite{Golub}.

Finally, we note that the numerical procedure developed here does not
apply to the $Q$ special eigenvalues of $A$, however as these are given
by Corollary \ref{cor:main-real} with exponential (in $n$) accuracy, this is not a
major limitation.

\vskip .2in\noindent
\begin{centering}
\section{Applications}
\label{chap:applications}
\end{centering}
\setcounter{figure}{0} \setcounter{equation}{0}

\noindent
Matrices like the one we study are often employed in systems of ordinary differential
equations. One example of this is in the study of traffic. If one assumes that the acceleration
of a car depends linearly on the perception of the relative velocities and positions of the
car in front of it and of the car behind it, then some analysis gives rise
to the equations
\begin{equation*}
\ddot x = B_1x+B_2\dot x,
\end{equation*}
where $B_1$ and $B_2$ are matrices of the type given in equation (\ref{eq:Bdef}) with the additional
property that they have row sum zero. In the special case that $B_1$ and $B_2$ are simultaneously
diagonalizable, one may use methods similar to those in this paper to study stability
(see for example \cite{tangerman}). In the more general case one takes refuge in the method
of \emph{periodic boundary conditions}. This raises the broader question of the
mathematical foundation of the validity of that method. Below we make some remarks
that relate that question to our present topic.

Discretizations of second order linear partial differential equations (PDE) naturally
give rise to tridiagonal matrices similar to the ones in this paper.
Below we give an example in 1 dimension. Here we notice that the theory in principle can
also be used in certain higher dimensional situations. Suppose we have a linear
second order PDE on a rectangle. Then we can discretize horizontally and vertically so that
each lattice point interacts with its horizontal neighbors through a matrix, say $L_1$,
and with its vertical neighbors through $L_2$. It is easy to show that the interaction
on the entire lattice is given by (see \cite{Golub}, section 4.8)
\begin{equation*}
L = L_1 \otimes I + I \otimes L_2
\end{equation*}
where $\otimes$ is the Kronecker product. The eigenvalues of $L$ are give by the Minkowski
sum of the eigenvalues of $L_1$ and $L_2$:
\begin{equation*}
\sigma(L)=\{\,z_1+z_2\,\big|\, z_1\in \sigma(A_1),\; z_2\in \sigma(A_2)\,\},
\end{equation*}
and the eigenvectors are given by the Kronecker product of the eigenvectors of $L_1$ and $L_2$.

We now make some more detailed comments.

\vskip .2in
\begin{centering}
\subsection*{Periodic Boundary Conditions}
\end{centering}

Perhaps the most common example of \emph{periodic boundary conditions} is part of the
foundation of solid state physics and has many applications in various technologies.
A set of identical ions on the line is separated by
a distance $a$ (a 1-dimensional \emph{Bravais} lattice). The position of the ion near
$ja$ is denoted by $x_j$ and is a function of time. After various approximations,
among which the assumption that ion interact only with their nearest neighbors, one
arrives at the following equation of motion:
\begin{equation}
\ddot x = q(A-2I)x
\label{eq:LinWithBdy}
\end{equation}
where $q$ is a positive constant related to the strength of the interaction and the mass of the ion.
Since physical systems are obviously finite, the remark is then, in the words of \cite{Ashcroft}
(Chapter 22), that
``we must specify how the ions at the two ends are to be described. [...] but this
would complicate the analysis without materially altering the final results. For if $N$
is large then [...] the precise way in which the ions at the ends are treated is
immaterial [...]". And thus one chooses a convenient way to
do that, namely periodic boundary conditions. The idea is clearly that \emph{physical
bulk} --- ie no boundary phenomena --- \emph{properties} are unchanged by the use of such boundary
conditions. That is: $A_{n+1}-2I$ is replaced by
\begin{equation*}
L_{n+1} =
\begin{pmatrix}
-2 & 1 & 0 & &\hdots &  1 \\
1 & -2 & 1 & &\hdots &   0 \\
\vdots &  & \ddots & \ddots& & \vdots\\
1 & \hdots & & 0 & 1 & -2
\end{pmatrix}\; ,
\end{equation*}
so that now:
\begin{equation}
\ddot x= qLx\;.
\label{eq:Laplacian}
\end{equation}
To the best of our knowledge, there is no mathematical proof for this
important fact at all.  It is thus tempting to employ the theory
developed here to have a closer look at this.

Write equations (\ref{eq:LinWithBdy}) and (\ref{eq:Laplacian}) as linear first order systems.
One easily sees that the eigenvalues of those systems are the roots of the eigenvalues
of $q(A-2I)$ and of $qL$, respectively. Suppose that the leading eigenvalues of the
systems are the most important ones for physical bulk properties. The eigenvalues of $L$
are well-known (and easily derived), namely $2(\cos\frac{2\pi m}{n+1}-1)$. Expanding this
to second order and taking the root, we obtain the (approximate) leading eigenvalues
for the system of equation (\ref{eq:Laplacian}):
\begin{equation*}
\nu_{m,\pm}= \pm i \,\frac{2 \pi m}{n+1} \quad m\in\{0,1,2,\cdots n\}\;.
\end{equation*}

We expect instabilities to fundamentally influence all physical properties. So if the matrix $A$ in
equation (\ref{eq:LinWithBdy}) has an eigenvalue with positive real part whereas the
matrix of $L$ in equation (\ref{eq:Laplacian}) does not, we can say that periodic boundary conditions
fails. That periodic boundary conditions \emph{can actually fail} in certain circumstances
is indicated by the following proposition.

\begin{figure}
\centerline{\includegraphics[width=5cm]{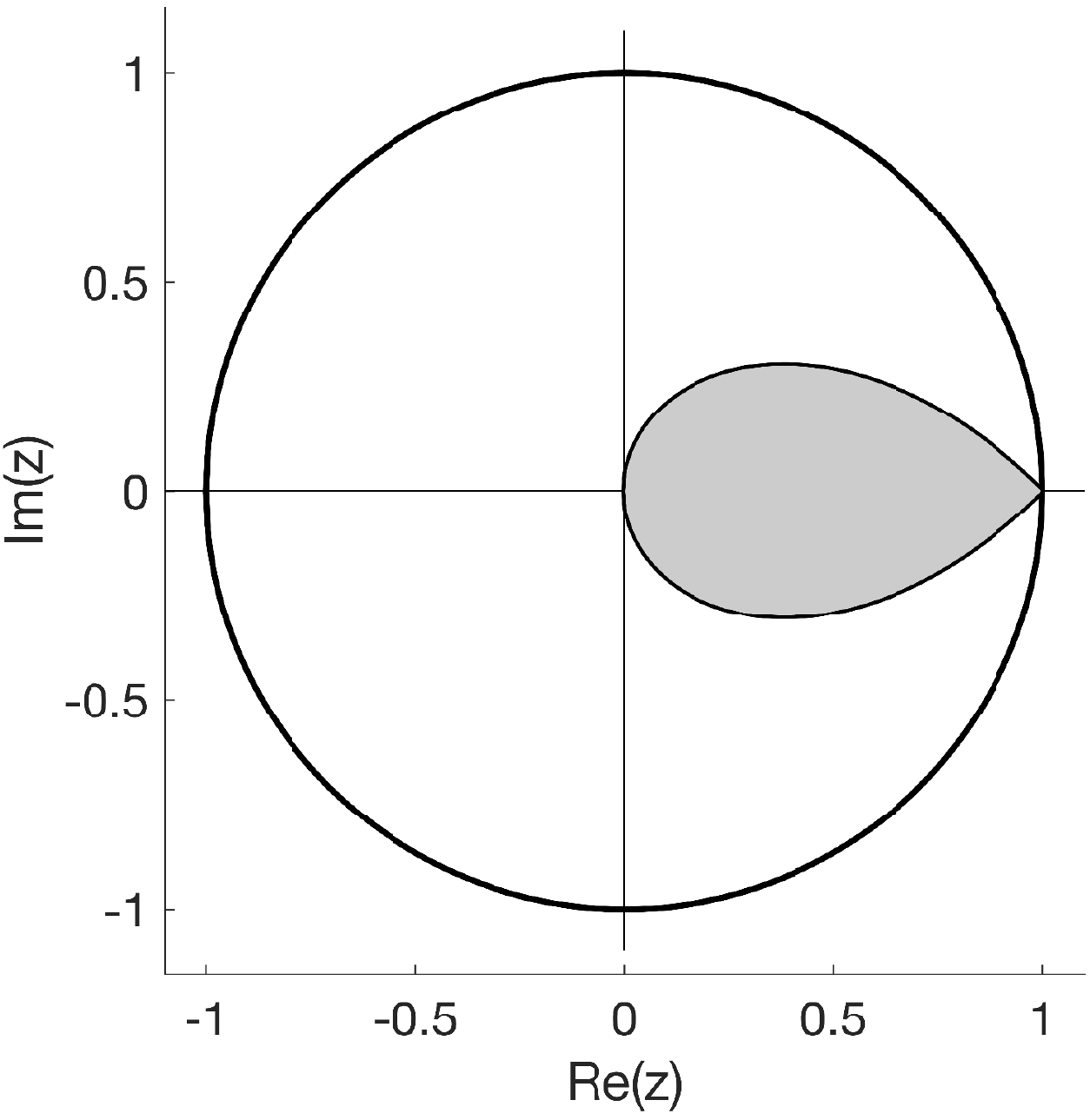}}
\caption{\emph{Shaded region consists of complex
  numbers $a+bi$ satisfying  $0<a<1$ and
  $|b|<\sqrt{\frac{a}{2-a}-a^2}$ from Theorem \ref{prop:failure}.}}
\label{fig:prop_failure}
\end{figure}

\vskip .2in\noindent
\begin{theo} Let $r_1=a+bi$ be such that $0<a<1$ and $|b|<\sqrt{\frac{a}{2-a}-a^2}$, and let
$r_2\in \mathbb{C}$ arbitrary (see figure \ref{fig:prop_failure}). If the matrix $A$ from equation (\ref{eq:Adef}) is such that
\begin{equation*}
b_0=-1/r_1 -1/r_2 \mbox{ and }b_1=\frac{1}{r_1 r_2}, \mbox{ with }c_0, c_{-1} \mbox{ arbitrary\;, }
\end{equation*}
and$\backslash$or
\begin{equation*}
c_{-1}= -1/r_1 -1/r_2 \mbox{ and  } c_0=\frac{1}{r_1 r_2}, \mbox{ with } b_0, b_1 \mbox{ arbitrary\;, }
\end{equation*}
then for sufficiently large $n$ the system corresponding to equation
(\ref{eq:LinWithBdy}) is unstable whereas the system with periodic
boundary conditions imposed (equation (\ref{eq:Laplacian})) is
Lyapunov stable.
\label{prop:failure}
\end{theo}

\noindent
{\bf Proof:} We have already shown that the eigenvalues of $L$ are
negative and thus its roots are imaginary (or 0) which makes the
system marginally stable. In fact, the dynamics is that of traveling waves (eg, see \cite{Ashcroft}),
which implies Lyapunov stability.

Assume the first condition holds. Simple algebra shows
$(\frac{z}{r_1}-1)(\frac{z}{r_2}-1)=b_1 z^2+ b_0 z + 1$, so $r_1$ and
$r_2$ are the roots of $b_1 z^2+ b_0 z + 1$, and are thus also roots
of $p(z)$ where $p(z)$ is given in Definition \ref{def:associated}. As
$r_1$ is inside the unit circle, Corollary \ref{cor:main} implies that
$A$ will have an eigenvalue exponentially close to $r_1+1/r_1$ (as
$n\rightarrow \infty$). Simple calculation shows
$Re(r_1+1/r_1) = a+\frac{a}{a^2+b^2}$ and that the inequality
$a+\frac{a}{a^2+b^2}>2$ is equivalent to the conditions on $r_1$ given
in the hypothesis. These conditions thus imply that the system
corresponding to equation (\ref{eq:LinWithBdy}) is unstable. A similar
argument applies if the second condition holds.
\QED

More detailed consideration of the notion of periodic boundary
conditions can be found in \cite{Cantos1}, \cite{Cantos2},
\cite{Herman}, and \cite{Herbrych}. However, it is still an open
question for what collection $C$ of possible boundary conditions the following holds.
For all boundary conditions in $C$, the leading eigenvalues of the Laplacian
are (close to) those of the system with periodic boundary conditions. Such
a statement would obviously of great value in all kinds of applications.
Even in the classical physics case, where the matrices are symmetric, as studied in \cite{Fonsecaea},
this problem is to the best of our knowledge unsolved (although the statement
is very widely used).

\begin{centering}
\subsection*{The Advection-Diffusion Equation}
\end{centering}

We consider a linear advection-diffusion equation on $[0,1]$
\begin{equation}
\partial_tu=\partial_x^2u + 2K\partial_x u\;,
\label{eq:advection-diffusion}
\end{equation}
with Dirichlet boundary conditions:
\begin{equation}
u(0,t)=f_0(t)\; , \quad u(1,t) = f_1(t)\; ,
\label{eq:dirichlet}
\end{equation}
and with Dirichlet-Neumann boundary conditions:
\begin{equation}
u(0,t)=f_0(t)\; , \quad \partial_x u(1,t) = f_1(t)\; .
\label{eq:dirichlet-neumann}
\end{equation}
Letting $u_j(t)$ stand for $u(\frac{j}{n},t)$ and using finite differences
(see \cite{Haberman}), one derives the following $n-1$-dimensional (not $n+1$ as in
the previous sections) system of ODE. Here $\dot u$ indicates derivative with respect to time of $u$. for the system with Dirichlet boundary conditions.
\begin{equation}
\begin{pmatrix}
\dot u_1\\
\dot u_2\\
\vdots\\
\vdots \\
\dot u_{n-1}\\
\end{pmatrix}
= n^2\,
\begin{pmatrix}
-2 & 1+\frac Kn & 0 & \hdots\\
1-\frac Kn & -2 & 1+\frac Kn & \hdots\\
\vdots & \ddots & & \ddots & \vdots\\
 & \hdots & 1-\frac Kn & -2 & 1+\frac Kn\\
 & \hdots & 0 & 1-\frac Kn& -2\\
\end{pmatrix}
\begin{pmatrix}
u_1\\
u_2\\
\vdots\\
\vdots \\
u_{n-1}\\
\end{pmatrix}
+ n^2\,
\begin{pmatrix}
(1-\frac Kn) f_0(t)\\
0\\
\vdots\\
0\\
(1-\frac Kn) f_1(t)\\
\end{pmatrix}\;.
\label{eq:discretization1}
\end{equation}
For the system with Dirichlet-Neumann boundary conditions, we obtain the following $n$ dimensional
system.
\begin{equation}
\begin{pmatrix}
\dot u_1\\
\dot u_2\\
\vdots\\
\vdots \\
\dot u_n\\
\end{pmatrix}
= n^2\,
\begin{pmatrix}
-2 & 1+\frac Kn & 0 & \hdots\\
1-\frac Kn & -2 & 1+\frac Kn & \hdots\\
\vdots & \ddots & & \ddots & \vdots\\
 & \hdots & 1-\frac Kn & -2 & 1+\frac Kn\\
 & \hdots & 0 & 2& -2\\
\end{pmatrix}
\begin{pmatrix}
u_1\\
u_2\\
\vdots\\
\vdots \\
u_n\\
\end{pmatrix}
+ n^2\,
\begin{pmatrix}
(1-\frac Kn) f_0(t)\\
0\\
\vdots\\
0\\
2n^{-1}(1+\frac Kn) f_1(t)\\
\end{pmatrix}\;.
\label{eq:discretization2}
\end{equation}
The matrix in this equation will be denoted by $B$. In the remainder of this section,
we are interested in the eigenvalues of the systems in equations (\ref{eq:discretization1})
and (\ref{eq:discretization2}).

\vskip .2in\noindent
\begin{prop} i: Fix $K$. Then for any $n>|K|$, all eigenvalues of the matrix $n^2B$ in
equation (\ref{eq:discretization1}) are real and less than $-K^2$.\\
ii: Fix $K\leq 0$. Then for any $n>|K|$, all eigenvalues of the matrix $n^2B$ in
equation (\ref{eq:discretization2}) are real and less than $-K^2$.
\label{prop:evalsPDE}
\end{prop}

\vskip .1in\noindent
{\bf Proof:} The proof of part i follows easily from that of part ii. We start with the latter.
First we use Appendix 1 to bring the matrix $B$ in the form used in this paper. Comparison with
equation (\ref{eq:Bdef}) shows that
\begin{equation*}
q\alpha_i^{-1} =1+\frac Kn \;,\;q\alpha_i=1-\frac Kn \;, \logand qd = -2\;.
\end{equation*}
Solve for $q$, $\alpha_i$ and $d$:
\begin{equation}
\alpha_i=\alpha \equiv \left(\frac{1-\frac Kn}{1+\frac Kn}\right)^{\frac12}\;,
q= \left(1-\frac{K^2}{n^2}\right)^{\frac12}\;, \logand
d= -2 \left(1-\frac{K^2}{n^2}\right)^{-\frac12}\;.
\label{eq:coeffsinPDE}
\end{equation}
Defining the diagonal matrix $D$ as in Lemma \ref{lem:app1-2}, one sees that
\begin{equation}
A\equiv D^{-1}\left(q^{-1}B-dI\right)D=
\begin{pmatrix}
0 & 1 & 0 & \hdots\\
1 & 0 & 1 & \hdots\\
\vdots & \ddots & & \ddots & \vdots\\
 & \hdots & 1 & 0 & 1\\
 & \hdots & 0 & \frac{2}{1-\frac Kn}& 0
\end{pmatrix}\;.
\label{eq:reduced-discretization}
\end{equation}
Comparison with equation (\ref{eq:Adef}) shows that
\begin{equation}
b_0=0 \; , \; b_1=0\; ,\; c_0=0\; , \logand c_{-1}=1 - \frac{2}{1-\frac Kn}=
-\frac{1+\frac Kn}{1-\frac Kn}=-\alpha^{-2}\;.
\label{eq:constantsinPDE}
\end{equation}
Thus the associated polynomial (see Definition \ref{def:associated}) is:
\begin{equation}
H(z) = z^{2n-2}\, z^2 \left(c_{-1}+ z^2\right) - \left( c_{-1} z^2 + 1\right)=
z^{2n} \left(z^2-\alpha^{-2}\right) - \left( 1-\alpha^{-2} z^2\right)\;.
\label{eq:HinPDE}
\end{equation}
Re-interpret the polynomial $p(z)=1-\alpha^{-2} z^2$ and the auxiliary
functions $f$ and $g$ as
\begin{equation}
f(z)=z^{2n} \quad \logand \quad g(z)=\frac{p(z)}{z^2p(z^{-1})}=
\frac{(1-\alpha^{-2} z^2)}{(z^2-\alpha^{-2})} \;.
\label{fandginPDE}
\end{equation}

We know from Proposition \ref{prop:phase_roots_real} that if $t_k$ is
a phase root, then $e^{it_k}$ is a root of $H(z)$. By Proposition \ref{prop:y2neq_deriv},
the roots $e^{it_k}$ of $H(z)/(z^2-1)$ correspond to eigenvalues
$\lambda_k=2\cos(t_k)$ of $A$.
By Corollary \ref{cor:app1}, the corresponding eigenvalues $\nu_k$ of
$n^2B$ are given by:
\begin{equation}
\nu_k = n^2(q\lambda_k+qd)=2n^2\left(1-\frac{K^2}{n^2}\right)^{\frac{1}{2}}\,\cos t_k - 2n^2<-K^2\;.
\label{eq:reg-evals-of-B}
\end{equation}
Therefore all eigenvalues of $n^2B$ arising in this way from roots of
$H(z)$ on the unit circle are real numbers less than zero, and no
instability arises from them.

When $K=0$, it follows that $\alpha=1$ and $H(z)$ simplifies to
\begin{equation*}
H(z)=(z^{2n}+1)(z^2-1)\;.
\end{equation*}
Clearly all of the roots of $H(z)/(z^2-1)$ lie on the unit circle and none are equal to 1,
and so all of the corresponding eigenvalues of $B$ are real numbers less than 0.

When $K<0$, then $\alpha>1$, and $p(z)$ (in equation
(\ref{fandginPDE})) has no roots \emph{inside} the unit circle, and so
$Q=0$. Adapting the proof of Lemma \ref{lem:windingnr} to our
re-interpreted $g(z)$ (by recognizing that $g(z)$ is now a rational
function of degree 2 rather than of degree 4) shows the winding number
of $g$ is $2Q-2$. A similar adaptation of Proposition
\ref{prop:phase-roots} shows that there must be at least $2n+2-2Q =
2n+2$ phase roots, each yielding a root of $H(z)$ on the unit
circle. But as $H(z)$ is a $2n+2$ degree polynomial, all of the roots
of $H(z)$ are on the unit circle and so all of the eigenvalues of $n^2 B$
are real numbers less than $-K^2$.

Now we return to part i. By the same reasoning as before,
we now obtain that $c_{-1}$ is also 0. Thus in this case,
\begin{equation*}
H(z) = z^{2n}-1 \;.
\end{equation*}
The roots of $H$ equal $e^{\frac{\pi i k}{n}}$, and are thus regular. The eigenvalues $\lambda_k$
of $A$ equal $\cos \frac{\pi k}{n}$  for $k\in\{1,\cdots n-1\}$ and the corresponding eigenvalues
of $n^2B$ are less than $-K^2$ as follows from equation (\ref{eq:reg-evals-of-B}). \QED

\vskip .1in
We return to the system with mixed boundary conditions. For $K>0$, $p(z)$ in equation (\ref{fandginPDE}) has precisely 2 roots inside
the unit circle. For symmetry reasons, these must be either on the unit circle, in which
case the corresponding eigenvalues of $n^2B$ are again less than $-K^2$, or else
the two roots are on the real line. In the latter case, one is the negative of the other.
This leads to two special eigenvalues $\lambda$ and $-\lambda$ of $A$, where $\lambda$ is a positive
real. By equation (\ref{eq:reg-evals-of-B}), we see that the eigenvalue of $n^2B$ which
corresponds to $-\lambda$ will tend to $-\infty$ as $n$ tends to $\infty$, and it is
therefore not relevant for the dynamics of the system.
One can show that the other special eigenvalue of $n^2B$ is always a real number in $(-K^2,0)$.
For brevity, we omit that argument.

\begin{figure}
\begin{center}
\begin{tabular}{cc}
\includegraphics[width=.45\linewidth]{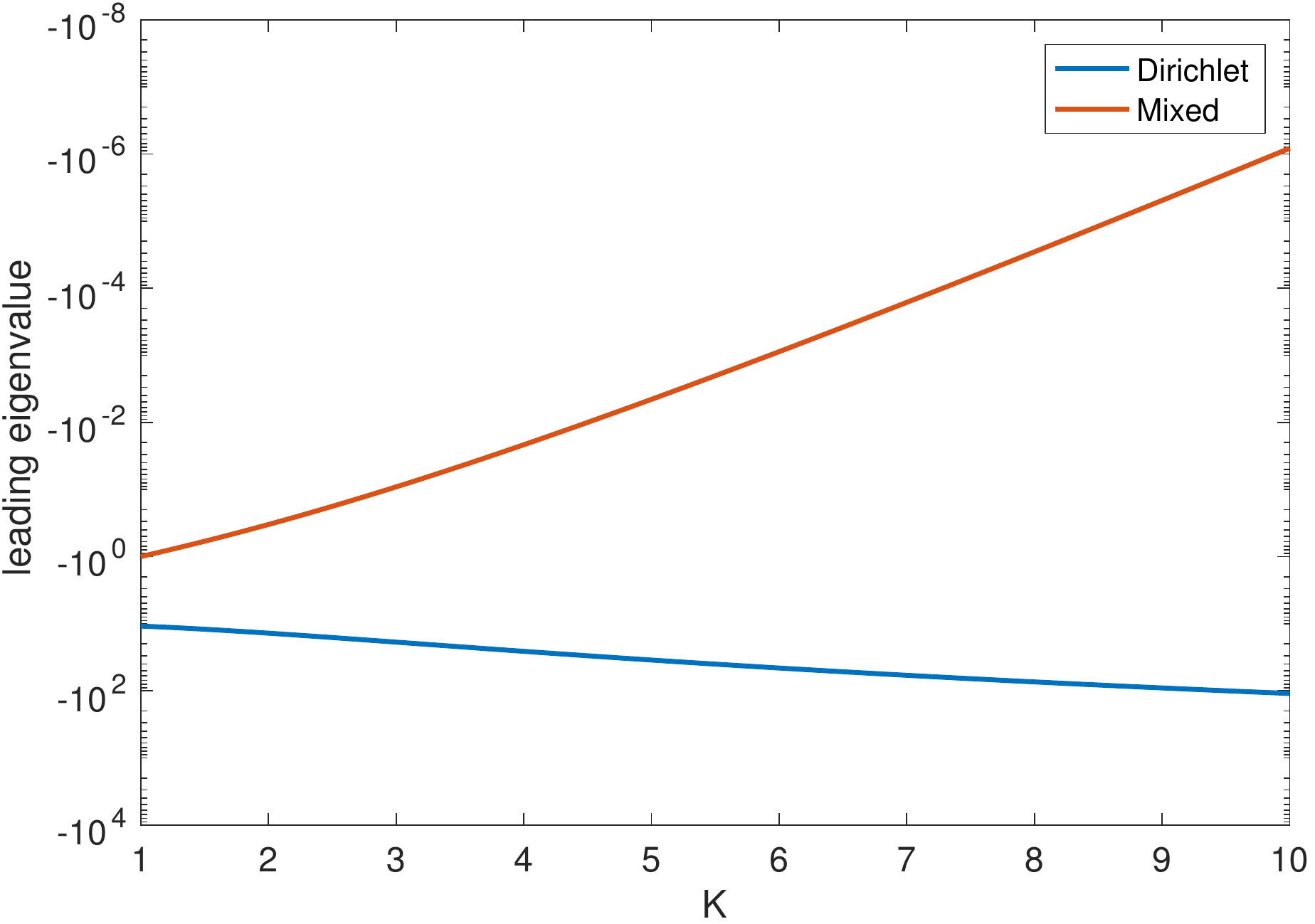} &
\includegraphics[width=.43\linewidth]{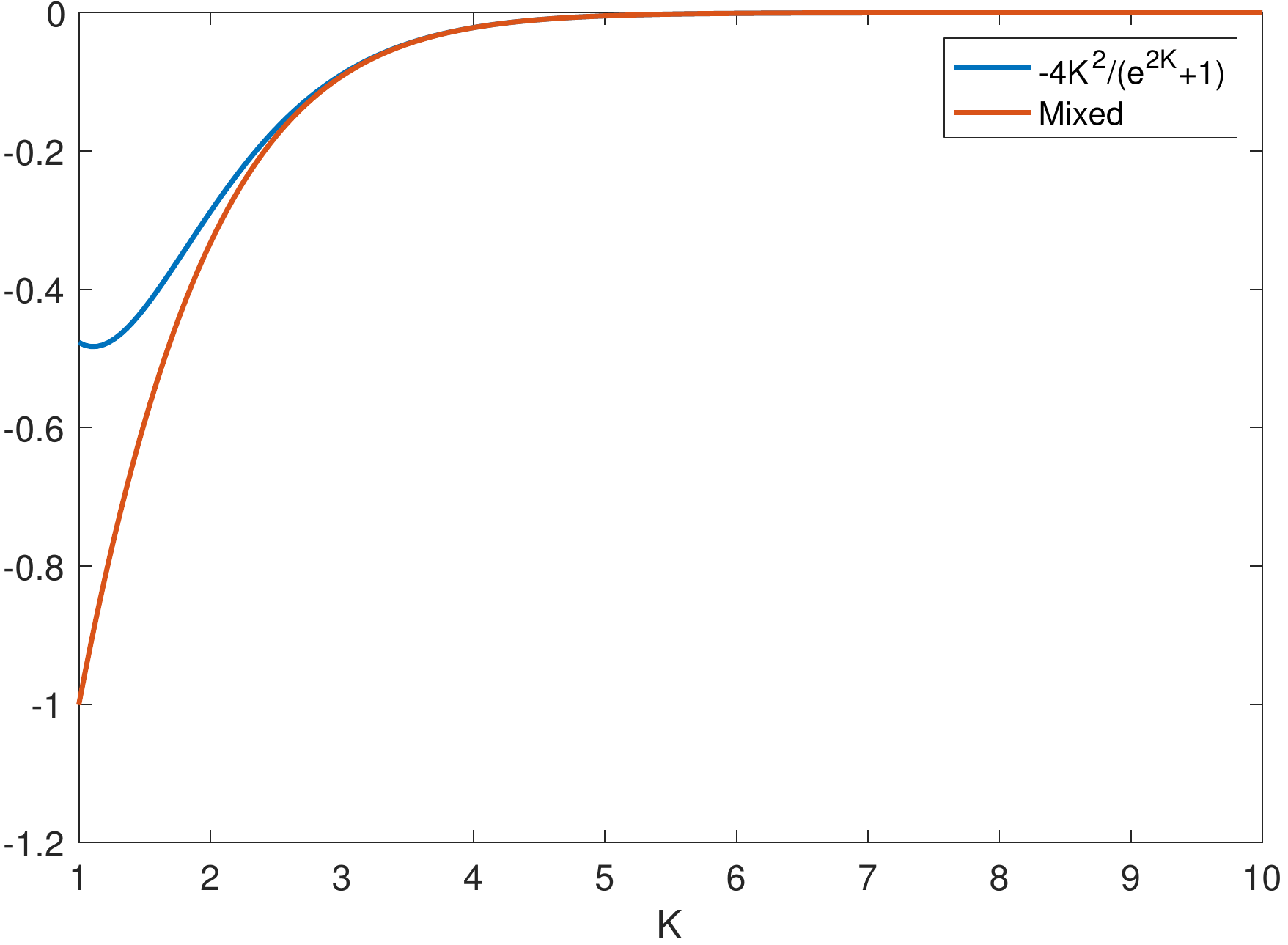} \\
(a) & (b) \\
\end{tabular}
\end{center}
\caption{ \emph{(a) Leading eigenvalues of $n^2 B$ for the mixed boundary
  conditions and the Dirichlet boundary conditions, vs K (note
  logarithmic scale) (b) Leading eigenvalue of $n^2 B$ for the mixed
  boundary condition system, and the prediction from Theorem
  \ref{theo:leading-eval}, convergence is observed as $K$ increases.}}
  \label{fig:leading-eval}
\end{figure}

Instead we will show that for large \emph{positive} $K$, the leading
eigenvalues of the two systems of equations (\ref{eq:discretization1})
and (\ref{eq:discretization1}) are very different. This is illustrated
in figure \ref{fig:leading-eval}. This implies a difference in
global dynamics (if given appropriate initial conditions) entirely due
to the different boundary conditions.

\vskip 0.2in\noindent
\begin{theo}
Let $K$ be positive and large. The leading eigenvalue of the system with
Dirichlet-Neumann boundary conditions is real and satisfies
\begin{equation*}
\nu= -\frac{4K^2}{e^{2K}+1} + {\cal O}\left(\frac{K}{\left( e^{2K}+1\right)^2}\right)
+ {\cal O}(n^{-2}) \;,
\end{equation*}
while the leading eigenvalue of the system with Dirichlet boundary conditions equals
$-K^2-\pi^2+ {\cal O}(n^{-2})$.
\label{theo:leading-eval}
\end{theo}

\vskip .0in\noindent
{\bf Proof:} The second part follows immediately from the previous proposition.

Fix a large value of $K$, we locate real roots of $H(z)$ for $n$ arbitrarily large.
To do this, set $\zeta\equiv z^{2n}$, $h(\zeta,s)\equiv H(z)$, and $s\equiv 1/n$, then
\begin{equation*}
h(\zeta,s)= \zeta \left(\zeta^s-\frac{1+ Ks}{1-Ks}\right)+\zeta^s \frac{1+Ks}{1-Ks}-1\;.
\end{equation*}
The equation for the corresponding eigenvalue of $n^2B$ becomes:
\begin{equation*}
\nu=s^{-2}\left(\sqrt{1-{K^2}{s^2}}\,\left(\zeta^{\frac s2}+\zeta^{-\frac s2}\right)-2\right)\;.
\end{equation*}
These equations have a meaningful expansions around $s=0$, namely
\begin{equation*}
h(\zeta,s)= \left[\zeta(\ln \zeta-2K)+(\ln \zeta+2K)\right]s+{\cal O}(s^2)\;.
\end{equation*}
Thus, in order for this equation to yield zero near $s=0$, we must have
\begin{equation}
\zeta(\ln \zeta-2K)+(\ln \zeta+2K)=0 \quad \logor \quad \ln \zeta (\zeta+1)=2K(\zeta-1)\;.
\label{eq:zeta}
\end{equation}
The expansion of the second equation (the eigenvalue) is
\begin{equation}
\nu = -K^2 + \left( \frac{\ln \zeta}{2}\right)^2+{\cal O}(s^2)\;.
\label{eq:eval-estimate}
\end{equation}

From equation (\ref{eq:zeta}) we see that if $K$ is positive and large and $\zeta \in(0,1)$, then
$\ln \zeta \approx -2K$, and thus $\zeta$ is very small. We make the following substitution
\begin{equation*}
u\equiv \frac 1K \quad \logand \quad \mu\equiv -u^{-2}+  \left( \frac{\ln \zeta}{2}\right)^2\;,
\end{equation*}
and obtain
\begin{equation*}
\ln \zeta = 2\sqrt{\mu+u^{-2}} \quad \logand \quad \zeta =  e^{2\sqrt{\mu+u^{-2}}}\;,
\end{equation*}
where $\mu$ is small when $u$ is small. Equation (\ref{eq:zeta}) becomes:
\begin{equation*}
 2\sqrt{\mu+u^{-2}}\left(e^{2\sqrt{\mu+u^{-2}}}+1\right) + 2u^{-1}
 \left( 1-e^{2\sqrt{\mu+u^{-2}}}\right) = 0 \;.
\end{equation*}
Multiplying by $\frac u2$ and rearranging gives:
\begin{equation*}
\sqrt{1+\mu u^2}=1-\frac{2}{e^{2\sqrt{\mu+u^{-2}}}+1}\;.
\end{equation*}
Note that as $K=u^{-1}$ becomes large, $\mu$ tends to zero exponentially in $K$. Squaring and then
subtracting 1, gives
\begin{equation}
\mu u^2 = -\frac{4}{e^{2 \sqrt{\mu+u^{-2}}+1}}+\frac{4}{\left(e^{2 \sqrt{\mu+u^{-2}}+1}\right)^2}\;.
\label{eq:to-be-inserted}
\end{equation}
Taylor expand the right hand side of this equation around $\mu=0$. Then substitute the first approximation for $\mu$. Finally, by equation (\ref{eq:eval-estimate}), $\nu$ and $\mu$
differ by ${\cal O}(n^{-2})$.

The proof of the second statement follows immediately from the proof of part i of Proposition
\ref{prop:evalsPDE}. Indeed, the reasoning there implies that $\nu_1 =
2n^2\left(1-\frac{K^2}{n^2}\right)^{\frac{1}{2}}\,\cos \frac{\pi}{n} - 2n^2$ which implies
the result.
\QED

\vskip .2in\noindent
\begin{centering}
\section{Appendix 1: More General Form of Matrices}
\label{chap:appendix1}
\end{centering}
\setcounter{figure}{0} \setcounter{equation}{0}

\noindent
In this appendix we show that with a little work one can expand the class of matrices to
which Corollary \ref{cor:main} can be applied. Namely, let $d$, $q$, and $\{\alpha_i\}_{i=1}^n$
be arbitrary complex numbers such that $k$ and $\alpha_i$ (for all $i$) are not zero. Define
\begin{equation}
{B}_{n+1} = q
\begin{pmatrix}
d-{b}_0 & \alpha_1^{-1}(1-{b}_1) & 0 & &\hdots &  0 \\
\alpha_1 & d & \alpha_2^{-1} & &\hdots &   0 \\
0    & \alpha_2 & d &  & \hdots & 0 \\
\vdots &  & \ddots & \ddots& & \vdots\\
0 & \hdots & & \alpha_{n-1} & d & \alpha_n^{-1} \\
0 & \hdots & & 0 & \alpha_n(1-{c}_{-1}) & d-{c}_0 \\
\end{pmatrix} .
\label{eq:Bdef}
\end{equation}
We characterize the eigenpairs of $B$ in terms of those of $A$ given in the main text.
For convenience of notation we drop the subscript $n+1$. The following results are simple computations.

\vskip .2in\noindent
\begin{lem} Let $D\in {\cal M}$ be the diagonal matrix with $\epsilon_i\neq 0$ as its $i$-th
diagonal element. Let $M\in {\cal M}$ arbitrary. Then
\begin{equation*}
\left(D^{-1}MD\right)_{ij} = \epsilon_i^{-1}\epsilon_j M_{ij} \;.
\end{equation*}
\label{lem:app1-1}
\end{lem}

\vskip -0.4in\noindent
\begin{lem} Set $\epsilon_i= \prod_{\ell<i}\,\alpha_\ell$ and $\epsilon_1=1$, and let $A$ be
the matrix given in equation (\ref{eq:Adef}). Then
\begin{equation*}
D^{-1}\left(q^{-1}B-dI\right)D=A \quad \logor \quad
B= q\left( DAD^{-1} + dI \right) \;.
\end{equation*}
\label{lem:app1-2}
\end{lem}

\vskip -0.4in\noindent
\begin{cory} Let $\lambda \in \mathbb{C}$ and $v\in \mathbb{C}^{n+1}$. Then
$\left(q(\lambda+d), Dv\right)$ is an eigenpair of $B$ if and only if $(\lambda,v)$ is an
eigenpair of $A$.
\label{cor:app1}
\end{cory}

\noindent
{\bf Remark:} In numerical work it is advantageous to work with the matrix $A$ and not with $B$,
because $B$ tends to have exponentially large condition number. This expresses itself in the
fact that regular eigenvectors $v$ of $A$ tend to have bounded components and, in contrast,
the regular eigenvectors $Dv$ of $B$ (see Lemma \ref{lem:app1-2}) tend to have components whose
ratios diverge as $\prod_{\ell<i}\,\alpha_\ell$. Clearly this can grow exponentially in $n$,
for example if all or most of the $\alpha_i>1+c$ and $c>0$.

We briefly mention two examples. Let $T\in {\cal M}$ be the tridiagonal Toeplitz matrix whose diagonal
elements equal $\delta$, whose sub-diagonal elements are equal to $\sigma$, and whose super-diagonal
elements are equal to to $\tau$. On the other hand, let $A_0\in {\cal
M}$ be the matrix
whose sub- and super-diagonal elements are 1, with $0$ on the diagonal. Corollary
\ref{cor:app1} says that
\begin{equation*}
T= \sqrt{\sigma\tau}\,DA_0D^{-1}+\delta I .
\end{equation*}
Since the spectrum of $A_0$ is easy to derive (namely,
$2\cos(\frac{\pi i}{n+2})$ for $i\in\{1,\cdots n+1\}$), the spectrum
of $T$ follows immediately (see \cite{noschese}).

For applications related to consensus forming and flocking, the following matrix was studied
in \cite{Hammond}:
\begin{equation}
{L}_{n+1} =
\begin{pmatrix}
\psi & 0 & 0 & &\hdots &  0 \\
\sigma & 0 & \tau & &\hdots &   0 \\
0    & \sigma & 0 &  & \hdots & 0 \\
\vdots &  & \ddots & \ddots& & \vdots\\
0 & \hdots & & \sigma & 0 & \tau \\
0 & \hdots & & 0 & \sigma+\phi & \theta\\
\end{pmatrix} \;.
\label{eq:Ldef1}
\end{equation}
One sees that $L$ is conjugate to $\tilde A$ where
\begin{equation}
{\tilde A}_{n+1} =
\sqrt{\sigma \tau}\,
\begin{pmatrix}
\frac{\psi}{\sqrt{\sigma \tau}} & 0 & 0 & &\hdots &  0 \\
1 & 0 & 1 & &\hdots &   0 \\
0    & 1 & 0 &  & \hdots & 0 \\
\vdots &  & \ddots & \ddots& & \vdots\\
0 & \hdots & & 1 & 0 & 1 \\
0 & \hdots & & 0 & 1-\frac\phi\sigma & \frac{\theta}{\sqrt{\sigma \tau}}\\
\end{pmatrix} \; .
\label{eq:Ldef2}
\end{equation}
The spectrum of $\tilde A$ can be studied with the methods of the main text.


\vspace{\fill}
\end{document}